\documentclass[a4paper,11pt]{amsart}
\usepackage{amssymb}
\usepackage{mathrsfs}
\usepackage{amsmath,amssymb,amsthm,latexsym,amscd,mathrsfs}
\usepackage{indentfirst}
 \newcommand{\ROM}[1]{\mathrm{\uppercase\expandafter{\romannumeral#1}}}
  \theoremstyle{definition}
 \newtheorem{thm}{Theorem}[section]

\newtheorem{conj}[thm]{Conjecture}


\title[A survey on the DDVV conjecture]{\textbf{A survey on the DDVV conjecture}}

\author[Ge J Q]{Ge Jianquan}\address{School of Mathematical Sciences, Laboratory of Mathematics and Complex Systems, Beijing Normal
University, Beijing 100875}\email{jqge@bnu.edu.cn}
\thanks {The project is partially supported by the NSFC ( No.10531090 and No.10229101 ) and the Chang Jiang Scholars
Program.}
\author[Tang Z Z]{Tang Zizhou}\address{School of Mathematical Sciences, Laboratory of Mathematics and Complex Systems, Beijing Normal
University, Beijing 100875}\email{zztang@mx.cei.gov.cn}

\thanks{The second author is the corresponding author.}
\date{}
\begin{document}
\maketitle

 \begin{center}
 Dedicated to Professor John C. Wood on his 60th birthday
 \end{center}

\section{Introduction}

Let $M^n$ be an immersed submanifold of a real space
form $N^{n+m}(c)$ of constant sectional curvature $c$. Let $R$
(resp. $\tilde{R}$) be the Riemann curvature tensor of $M^n$ (resp.
$N^{n+m}(c)$), $h$ the second fundamental form, $A_{\xi}$ the
shape-operator associated to a normal vector field $\xi$, and
$R^{\bot}$ the curvature tensor of the normal connection. The
equations of Gauss and Ricci are given by
\[\langle R(X,Y)Z, T\rangle=\langle \tilde{R}(X,Y)Z, T\rangle+\langle h(X,T),h(Y,Z)\rangle-\langle h(X,Z), h(Y,T)\rangle,\]
\[\langle R^{\bot}(X,Y)\xi, \eta\rangle=\langle \tilde{R}(X,Y)\xi, \eta\rangle+\langle [A_{\xi},A_{\eta}]X, Y\rangle,\]
for tangent vectors $X$,$Y$,$Z$,$T$ and normal vectors
$\xi$,$\eta$.

Let $\{e_1,...,e_n\}$ (resp. $\{\xi_1,...,\xi_m\}$)
be an orthonormal basis of $T_pM$ (resp. $T^{\bot}_pM$). The
normalized scalar curvature $\rho$ and normal scalar curvature
$\rho^{\bot}$ of $M^n$ at p are defined by
\[\rho=\frac{2}{n(n-1)}\sum^n_{1=i<j}\langle R(e_i,e_j)e_j,e_i)\rangle,\]
\[\rho^{\perp}=\frac{2}{n(n-1)}\left(\sum^n_{1=i<j}\sum^m_{1=r<s}\langle R^{\perp}(e_i,e_j)\xi_r, \xi_s\rangle^2\right)^{\frac{1}{2}}=\frac{2}{n(n-1)}|R^{\bot}|.\]

 Let
$H=\frac{1}{n}Tr(h)=\frac{1}{n}\sum_{i=1}^nh(e_i,e_i)$ be the mean
curvature vector field.
\begin{thm}(\cite{W}, \cite{GR})
Let $M^2$ be an immersed surface of a real space form $N^{2+m}(c)$.
Then \[\rho+\rho^{\bot}\leq |H|^2+c\] at every point p of $M^2$,
with equality if and only if the ellipse of curvature is a circle.
\end{thm}
This inequality was proved by Wintgen\cite{W} in 1979 for surfaces
in $4$ dimensional Euclidean space and by Guadelupe and
Rodriguez\cite{GR} in 1983 for surfaces in arbitrary space form. By
applying this inequality, they derived many rigidity results about
immersions of surfaces including some results proved by
Barbosa\cite{B}, Yau\cite{Y}, \cite{Y2}, Ruh\cite{R}. Since the
equality in the inequality holds if and only if the ellipse of
curvature is a circle, and the property that the ellipse is a circle
is a conformal invariant, surfaces in $\mathbb{R}^4$ that attain
equality everywhere could be stereographically projected from those
in $S^4$ (for instance see do Carmo and Wallach\cite{CW}). In
\cite{GR}, they also remarked that Atiyah and Lawson (unpublished)
had shown that an immersed surface in $S^4$ has the ellipse always a
circle if and only if the canonical lift of the immersion map into
the bundle of almost complex structures of $S^4$ (\emph{i.e.},
$CP^3$) is holomorphic (also see Theorem 5 in \cite{BFLPP}). One can
also consider surfaces with circular ellipse of curvature in general
Riemannian manifolds other than real space forms, in which case
circular ellipse of curvature also corresponds to the equality
condition of some similar inequality involving relative curvatures.
For this subject on surfaces with circular ellipse of curvature, we
refer to \cite{BFLPP}, \cite{Dajczer-T1}, \cite{Castro} and
references therein.

In 1996, Chen\cite{ChenBY0} proved the following similar but more
elementary inequality for submanifolds in real space forms (also see
Suceav$\breve{a}$\cite{Suceava}, Ge\cite{Ge1}):
$$\rho\leq|H|^2+c.$$

To study such inequalities relating intrinsic and extrinsic
curvature invariants, in 1999, De Smet, Dillen, Verstraelen and
Vrancken\cite{DDVV} defined the normal scalar curvature for high
dimensional submanifolds and posed a conjecture as generalization of
the previous Theorem (which is now called the DDVV conjecture):
\begin{conj}\label{conj1}(\cite{DDVV})
Let $M^n$ be an immersed submanifold of a real space form
$N^{n+m}(c)$. Then \[\rho+\rho^{\bot}\leq |H|^2+c.\]
\end{conj}
The conjecture was proved for $m=2$ by them, where also some partial
classification results were obtained in case equality holds at every
point.

 Recent developments mostly based on the following
translation of the conjecture established by Dillen, Fastenakels and
Veken\cite{DFV}.
\begin{thm}\label{conj2}(\cite{DFV})
Conjecture 1.2 is true for submanifolds of dimension $n$ and
codimension $m$ if for every set $\{B_1, . . . ,B_m\}$ of real
symmetric $(n\times n)$-matrices with trace zero the following
inequality holds: \[\sum_{r,s=1}^m\|[B_r,B_s]\|^2\leq
\left(\sum_{r=1}^m\|B_r\|^2\right)^2.\]
\end{thm}
In fact, putting $B_r=A_{\xi_r}-\langle H, \xi_r\rangle id$, and
using Gauss and Ricci equations, we find
\[|H|^2-\rho+c=\frac{1}{n(n-1)}\sum_{r=1}^m\|B_r\|^2,\]
\[\rho^{\bot}=\frac{1}{n(n-1)}\left(\sum_{r,s=1}^m\|[B_r,B_s]\|^2\right)^{\frac{1}{2}}.\]

When $m=2$, such kind of matrix inequality in Theorem \ref{conj2}
was already studied by Simons\cite{Simons} and Chern\cite{chern} in
proving the well-known scalar curvature pinching theorem for minimal
submanifolds in spheres.

\section{Recent developments}
For convenience, we denote both inequalities in Conjecture
\ref{conj1} and Theorem \ref{conj2} by $P(n,m)$. Thus, $P(2,m)$ and
$P(n,2)$ were proved in \cite{W}, \cite{GR} and \cite{DDVV}
respectively as mentioned before. We'd like to list the recent
developments as follows:
 \begin{itemize}
    \item  In 2004, Dillen, Haesen, Petrovi$\acute{c}$-Torga$\check{s}$ev
and Verstraelen\cite{DHTV} proved the $P(n,2)$ of Lorentz type,
\emph{i.e.}, for Lorentzian manifold embedded locally and
isometrically in a pseudo-Euclidean space with codimension $2$.

    \item  A weaker version of the conjecture and also a special case when the
    submanifold is either H-umbilical Lagrangian in $\mathbb{C}^n$ or ultra-minimal in
    $\mathbb{C}^4$ were proved in the same paper as Theorem \ref{conj2} by Dillen, Fastenakels and
Veken\cite{DFV}.

    \item  The first nontrivial case $P(3,m)$ was proved by Choi and Lu\cite{CL} in 2006, where
    they also classified all $3$ dimensional minimal submanifolds satisfying the equality everywhere
    and compared this conjecture with the comass problem in calibrated geometry.

    \item  Another nontrivial case $P(n, 3)$ was proved by Lu\cite{Lu1} in 2007, where he also discussed the relationship between DDVV conjecture and pinching theorems of scalar curvature for minimal submanifolds in spheres.

    \item  Besides the introduction of some recent developments,  Lu\cite{Lu2} discussed some applications of DDVV conjecture and its relation
    with a conjecture of B$\ddot{o}$ttcher and Wenzel\cite{Bottcher-W} in random matrix theory.

    \item  In 2007, the whole conjecture was proved by Lu\cite{Lu3} and the authors\cite{GT1}
     independently by quite different methods. The equality condition (locally) was also determined completely by
     us (see later).

    \item Lu\cite{Lu3} also observed some interesting applications of this
     conjecture to the Simons type pinching theorems (see Simons\cite{Simons}) and established
     a new pinching theorem, in view of which he made a conjecture about the gap theorem of Peng-Terng type (see Peng and Terng\cite{Peng-T}, Chern, do Carmo and Kobayashi\cite{CCK}) in high
     codimensions. There's also a proof of the conjecture of B$\ddot{o}$ttcher and Wenzel\cite{Bottcher-W} in \cite{Lu3}.

    \item By the same method (but a little more complicated) as \cite{GT1}, Ge\cite{Ge2} obtained a similar optimal inequality like $P(n,m)$ for real skew-symmetric matrices (see later).

     \item Some partial results on classification of submanifolds when the equality in the DDVV inequality holds everywhere
     were proved by Dajczer and Florit\cite{Dajczer-F}, Dajczer and Tojeiro\cite{Dajczer-T1}, \cite{Dajczer-T2}. In particular, \cite{Dajczer-T2} provided a parametric construction in terms of minimal surfaces of the
Euclidean submanifolds of codimension two and arbitrary dimension
that attain equality everywhere in the DDVV inequality. As they
remarked that by the pointwise equality condition determined in
\cite{GT1}, their results provide a complete classification of all
non-minimal submanifolds (of arbitrary codimension) of dimension at
least four that attain equality in the DDVV inequality and whose
first normal spaces have dimension two everywhere . It is a very
interesting problem to study the remaining cases.
\end{itemize}

We now state the main theorems in \cite{GT1} and \cite{Ge2} as follows.
\begin{thm}(\cite{Lu3},\cite{GT1})\label{thm1}
Let $B_1,...,B_m$ be $(n\times n)$ real symmetric matrices. Then
\[\sum_{r,s=1}^m\|[B_r,B_s]\|^2\leq (\sum_{r=1}^m\|B_r\|^2)^2,\]
where the equality holds if and only if under some $O(n)\times O(m)$
action\footnote{A $(P,R)\in O(n)\times O(m)$ acts on
$(B_1,\cdots,B_m)$ by: $(P,R)\cdot
(B_1,\cdots,B_m):=(P^tB_1P,\cdots,P^tB_mP)R$.} all matrices $B_r$
are zero except the following $2$ matrices:
\[H_1=\left(\begin{array}{ccccc}\mu&
0& 0&\cdots& 0\\0& -\mu& 0&\cdots& 0\\0& 0& 0&\cdots& 0\\
\vdots&\vdots & \vdots&\ddots &\vdots
\\0& 0& 0&\cdots& 0 \end{array}\right),\hskip 0.5cm H_2=\left(\begin{array}{ccccc}0& \mu& 0&\cdots& 0\\\mu& 0& 0&\cdots&
0\\0& 0& 0&\cdots& 0\\ \vdots&\vdots & \vdots&\ddots &\vdots
\\0& 0& 0&\cdots& 0 \end{array}\right),\]
where $\mu\geq 0$ is a real number.
\end{thm}
\begin{thm}(\cite{Lu3},\cite{GT1})
Let $f: M^n\rightarrow N^{n+m}(c)$ be an isometric immersion. Then
\[\rho+\rho^{\perp}\leq |H|^2+c,\]
where the equality holds at some point $p\in M$ if and only if there
exist an orthonormal basis $\{e_1,...,e_n\}$ of $T_pM$ and an
orthonormal basis $\{\xi_1,...,\xi_m\}$ of $T_p^{\perp}M$, such that
\[A_{\xi_1}=\left(\begin{array}{ccccc}\lambda_1+\mu& 0& 0&\cdots& 0\\0& \lambda_1-\mu& 0&\cdots&
0\\0& 0& \lambda_1&\cdots& 0\\ \vdots&\vdots & \vdots&\ddots &\vdots
\\0& 0& 0&\cdots& \lambda_1 \end{array}\right),\hskip 0.5cm A_{\xi_2}=\left(\begin{array}{ccccc}\lambda_2 & \mu& 0&\cdots&
0\\\mu& \lambda_2& 0&\cdots& 0\\0& 0& \lambda_2&\cdots& 0\\
\vdots&\vdots & \vdots&\ddots &\vdots
\\0& 0& 0&\cdots& \lambda_2 \end{array}\right),\]
$A_{\xi_3}=\lambda_3I_n$ and all other shape operators
$A_{\xi_r}=0$, where $\mu,\lambda_1,\lambda_2,\lambda_3$ are real
numbers.
\end{thm}
\begin{thm}(\cite{Ge2})
Let $B_1,...,B_m$ be $(n\times n)$ real skew-symmetric matrices. Then
\begin{itemize}
\item[(i)] for $n=3$,
 \[\sum_{r,s=1}^m\|[B_r,B_s]\|^2\leq
\frac{1}{3}\left(\sum_{r=1}^m\|B_r\|^2\right)^2,\] where the
equality holds if and only if under some $O(n)\times O(m)$ action
all matrices $B_r$ are zero except the following $3$ matrices:
\[C_1:=\left(\begin{array}{ccc}0&
\lambda& 0\\-\lambda& 0& 0\\0& 0& 0
\end{array}\right),~~C_2:=\left(\begin{array}{ccc}0&
0& \lambda\\0& 0& 0\\-\lambda& 0& 0
\end{array}\right),~~C_3:=\left(\begin{array}{ccc}0& 0& 0\\0&
0& \lambda\\0& -\lambda& 0
\end{array}\right),\]
where $\lambda\geq 0$ is a real number;

\item[(ii)] for $n\geq4$,
 \[\sum_{r,s=1}^m\|[B_r,B_s]\|^2\leq
\frac{2}{3}\left(\sum_{r=1}^m\|B_r\|^2\right)^2,\] where the
equality holds if and only if under some $O(n)\times O(m)$ action
all matrices $B_r$ are zero except the following $3$ matrices:
$diag(D_1, 0)$, $diag(D_2, 0)$, $diag(D_3, 0)$, where $0\in M(n-4)$
is a zero matrix, $\lambda\geq 0$ is a real number,
\[D_1:=\left(\begin{array}{cccc}0& \lambda& 0&0\\-\lambda& 0& 0&0\\0&
0& 0&\lambda\\0&0&-\lambda&0
\end{array}\right),
D_2:=\left(\begin{array}{cccc}0& 0&\lambda& 0\\0& 0&
0&-\lambda\\-\lambda& 0& 0&0\\0&\lambda&0&0
\end{array}\right),
D_3:=\left(\begin{array}{cccc}0& 0&0&\lambda\\0& 0&\lambda&0\\0&
-\lambda& 0&0\\-\lambda&0&0&0
\end{array}\right).\]
\end{itemize}
\end{thm}

\section{Sketch of our method}
Due to the success of \cite{Ge2} in finding and proving optimal
inequalities involving norms of commutators or Lie brackets by the
method of \cite{GT1}, we'd like to give a sketch of the proof of
Theorem \ref{thm1} to conclude this survey. It goes in $3$ major
steps as follows:
\begin{itemize}
  \item [$(1)$] {Translation of the inequality into a polynomial
  inequality}
\end{itemize}
Denote by $SM(n)$ the $N:=\frac{n(n+1)}{2}$-dimensional space of
real symmetric $(n\times n)$ matrices. Let $\{\hat{E}_{ij}\}_{1\leq
i\leq j\leq n}$ be the standard orthonormal basis of $SM(n)$. Taking
the letter order for the indices set $S:=\{(i,j)|1\leq i\leq j\leq
n\}$, we identify $S$ with $\{1,...,N\}$. Then there's a unique
$(N\times m)$ matrix $B$ such that
\[(B_1,...,B_m)=(\hat{E}_1,...,\hat{E}_N)B.\]
Similarly, an orthogonal matrix $Q\in SO(N)$ determines $N$ real
symmetric $(n\times n)$ matrices $\{\hat{Q}_1,...,\hat{Q}_N\}$ by
\[(\hat{Q}_1,...,\hat{Q}_N)=(\hat{E}_1,...,\hat{E}_N)Q.\]

Since $BB^t$ is a semi-positive definite matrix in $SM(N)$, there
exists an orthogonal matrix $Q\in SO(N)$ such that $BB^t=Q\hskip
0.1cm diag(x_1,...,x_N)\hskip 0.1cm Q^t$ with $x_{\alpha}\geq 0$,
$1\leq\alpha\leq N$. Then direct calculations show that
\[\sum_{r=1}^m\|B_r\|^2=\|B\|^2 =\sum_{\alpha=1}^Nx_{\alpha},\]
\[\sum_{r,s=1}^m\|[B_r, B_s]\|^2=\sum_{\alpha,\beta=1}^Nx_{\alpha}x_{\beta}\|[\hat{Q}_{\alpha},
\hat{Q}_{\beta}]\|^2.\] Now the problem is changed into proving the
following \underline{polynomial inequality}:
\[f_Q(x):=F(x,Q):=\sum_{\alpha,\beta=1}^Nx_{\alpha}x_{\beta}\|[\hat{Q}_{\alpha},
\hat{Q}_{\beta}]\|^2-\left(\sum_{\alpha=1}^Nx_{\alpha}\right)^2\leq0,\] for all
$x\in \mathbb{R}^N_{+}:=\{0\neq x=(x_1,...,x_N)\in
\mathbb{R}^N|x_{\alpha}\geq 0, 1\leq\alpha\leq N\}$ and $Q\in
SO(N)$.

 \begin{itemize}
  \item [$(2)$] {Approximate process}
\end{itemize}
Since $f_Q(x)$ is a homogeneous polynomial, we need only to consider
the inequality on the subset
$\bigtriangleup:=\{x\in\mathbb{R}^N_{+}|\sum_{\alpha}x_{\alpha}=1\}$
and show that \[G:=\{Q\in SO(N)|f_Q(x)\leq 0, \forall x\in
\bigtriangleup\}=SO(N).\] By an observation to the case when
$Q=I_N$, we find it may happen to be that $f_Q(x)<0$ for all $x$ in
the interior of $\bigtriangleup$, which inspires us to consider the
following approximate process. \\
For any sufficiently small $\varepsilon>0$, put
$\bigtriangleup_{\varepsilon}:=\{x\in \bigtriangleup|x_{\alpha}\geq
\varepsilon, 1\leq\alpha\leq N\}$ and $G_{\varepsilon}:=\{Q\in
SO(N)|f_Q(x)< 0, \forall x\in \bigtriangleup_{\varepsilon}\}$. Try
to prove \[G_{\varepsilon}=SO(N).\] Finally,
\[G=\lim_{\varepsilon\rightarrow 0}G_{\varepsilon}=SO(N).\]

\begin{itemize}
   \item [$(3)$] {Continuity method and a priori
   estimate}
\end{itemize}
Since $SO(N)$ is connected, we can use continuity method which
consists of the following three steps:
\\
\underline{step 1}: \hskip 0.3cm $I_N\in G_{\varepsilon}$,  (thus $G_{\varepsilon}\neq\emptyset$),\\
\underline{step 2}: \hskip 0.3cm $G_{\varepsilon}$ is open in $SO(N)$, \\
\underline{step 3}: \hskip 0.3cm $G_{\varepsilon}$ is closed in
$SO(N)$.\vskip 0.3cm
 Step 1 and Step 2 can be easily verified, while step 3 needs the
 following a priori estimate which is the main difficulty in the
 whole proof and indicates the idea of the proof of the equality
 case.\vskip 0.3cm
{\bf{a priori estimate}}: Suppose $f_Q(x)\leq 0$, for every $
x\in\bigtriangleup_{\varepsilon}$. Then $f_Q(x)< 0$, for every $
x\in\bigtriangleup_{\varepsilon}$.

The proof of the a priori estimate depends on a sequence of lemmas
mostly depending on the concrete properties of those coefficients of
Lie brackets of a given basis of the space in problem. Once one
knows more about those coefficients (structural constants when the
space is a Lie algebra), one will probably find and prove such type
optimal inequality by our method.

\end{document}